\documentclass[a4paper]{article}

\title{An infinitude of proofs of the infinitude of primes}

\author{L. J. P. Kilford}

\begin{document}

\maketitle
\begin{center}
\emph{Dedicated to Yang-Hui He, for his $30^{\rm th}$ birthday.}
\end{center}
Since Euclid's proof that there are infinitely many primes, there have been  many different proofs of this theorem; some using topology, others combinatorics or number theory. In this note, we present a generalization of a method of Perott (\cite{perott}, given in~\cite{ribenboim}, which references~\cite{dickson}, page~413) which gives infinitely many proofs of this theorem.

Let~$m \ge 2$ be an integer. It is well-known that the sum~$\sum_{n=1}^\infty (1/n^m)$ converges; moreover, it converges to a real number which is between~1 and~2, and is strictly less than~2, because~$\sum_{n=1}^\infty (1/n^2)=\frac{\pi^2}{6} < 2$.

We suppose that there are finitely many primes~$p_1 < p_2 <\cdots < p_r$, and let~$N$ be a positive integer which is larger than the product of all of these primes. We will now divide the natural numbers~$t \le N$ into two classes; those which are $m^{\rm th}$ power free, and those which are divisible by an $m^{\rm th}$ power. There are~$m^r$ $m^{\rm th}$ power free numbers less than or equal to~$N$, because this is the number of products of primes which are $m^{\rm th}$ power free. On the other hand, the number of integers~$t \le N$ which are divisible by~$p_i^m$ is at most~$N/p_i^m$, so the number of natural numbers~$t \le N$ divisible by an $m^{\rm th}$ power is at most~$\sum_{i=1}^r (N/p_i^m)$. Hence we have
\[
N \le m^r  + \sum_{i=1}^r \frac{N}{p_i^m} < m^r + N\left(\sum_{n=1}^\infty \frac{1}{n^m} - 1\right) = m^r + N(1-\delta),
\]
where~$\delta > 0$. We now choose~$N$ such that~$N \delta \ge m^r$, which we can do because~$m$ and~$r$ are fixed, and we obtain a contradiction. The infinitude of proofs follows because every natural number~$m \ge 2$ suffices here; Perrot's proof used~$m=2$.


\end{document}